\documentclass[reqno]{amsart}

\usepackage{graphicx,amssymb,mathrsfs,amsmath,color,fancyhdr}
\usepackage[all]{xy}
\newtheorem{theorem}{Theorem}[section]
\newtheorem{lemma}[theorem]{Lemma}

\theoremstyle{definition}
\newtheorem{definition}[theorem]{Definition}

\newtheorem{proposition}[theorem]{Proposition}

\theoremstyle{remark}
\newtheorem{remark}[theorem]{Remark}

\numberwithin{equation}{section}

\begin{document}

\title[New constructions of NMDS self-dual codes]
{New constructions of NMDS self-dual codes}

\author{Dongchun Han}
\address{Department of Mathematics, Southwest Jiaotong University, Chengdu 610000, P.R. China}
\email{handongchun@swjtu.edu.cn}
\author{Hanbin Zhang}
\address{School of Mathematics (Zhuhai), Sun Yat-sen University, Zhuhai 519082, Guangdong, P.R. China}
\email{zhanghb68@mail.sysu.edu.cn}


\begin{abstract}
Near maximum distance separable (NMDS) codes are important in finite geometry and coding theory. Self-dual codes are closely related to combinatorics, lattice theory, and have important application in cryptography.
In this paper, we construct a class of $q$-ary linear codes and prove that they are either MDS or NMDS which depends on certain zero-sum condition. In the NMDS case, we provide an effective approach to construct NMDS self-dual codes which largely extend known parameters of such codes. In particular, we proved that for square $q$, almost $q/8$ NMDS self-dual $q$-ary codes can be constructed.

\end{abstract}

\maketitle{}

\section{Introduction}

Let $\mathbb F_q$ be a finite field with $q$ elements of characteristic $p$.
A $q$-ary $[n,k,d]$-linear code $\mathcal C$ is a subspace of $\mathbb F_q^n$ with dimension $k$ and minimum distance $d$. For any $a=(a_1,\ldots,a_n)\in \mathbb F_q^n$ and $b=(b_1,\ldots,b_n)\in \mathbb F_q^n$, their Euclidean inner product is defined as $\langle a,b\rangle=\sum_{i=1}^na_ib_i$.
We call $\mathcal C^{\bot}=\{b\in\mathbb F_q^n\ |\ \langle a,b\rangle=0\text{ for all }a\in\mathcal C\}$ the dual code of $\mathcal C$. The famous Singleton bound states that $d\le n-k+1$. A $q$-ary $[n,k,d]$-linear code $\mathcal C$ is called  maximum distance separable (MDS) if $d=n-k+1$. The study of MDS codes is one of the central topics in coding theory. A $q$-ary $[n,k,d]$-linear code $\mathcal C$ is called  almost maximum distance separable (AMDS)  if $d=n-k$. An AMDS code $\mathcal C$ is called near maximum distance separable (NMDS) if $C^{\bot}$ is also AMDS; see \cite{DodL}. NMDS codes are closely related to design theory and finite geometry (see, e.g., \cite{DingT,DodL,WangH}), and also have application in secret sharing scheme \cite{SimV,ZhouWXLQY}.

A linear code $\mathcal C$ is called (Euclidean) self-dual if $\mathcal C=\mathcal C^{\bot}$. Self-dual codes have various connections with combinatorics and lattice theory; see, e.g., \cite{ConS,MacWS}. In practice, self-dual codes have important application in cryptography; see \cite{CraDGULMP,DouMS}.

Due to the aforementioned reasons, it is natural and important to study MDS self-dual codes. In fact, MDS self-dual codes have been extensively studied recently; see, e.g., \cite{FLL,FLLL,FangF,FangXF,JinX,KimLee,LLL,HYan,ZhangFeng}. Compared with the studies on MDS self-dual codes, the constructions on NMDS self-dual codes are not too much. Note that, the parameter $[n,k,d]$ of an NMDS self-dual code $\mathcal C$ is completely determined by its length $n$. Paralell to the famous MDS main conjecture (which states that the maximal length for $q$-ary nontrivial MDS codes is $q + 1$ except when $q$ is even and $k\in\{3, q-1\}$), it is conjectured in \cite{LanR} that for all prime powers $q$, the maximal length for $q$-ary nontrivial NMDS codes is at most $2q + 2$, which is twice of the MDS case.
In \cite{GKL}, $q$-ary NMDS self-dual codes of length $q-1$ for odd prime power $q$ with $q\equiv 1\pmod 4$ and $q\le 113$ were constructed via Reed-Solomon codes. In \cite{KotKS}, the existence of $p$-ary NMDS self-dual codes of small length (up to 16) over some small prime fields with $p\le 197$ was proved. Recently, using elliptic curves, Jin and Kan \cite{JinK} constructed more NMDS self-dual codes: let $q$ be an odd prime power satisfying $q\equiv 1\pmod 4$, then for every $n$ divided by 4 that satisfies $q\ge 4^{n+3}(n+3)^2$, there exists a $q$-ary NMDS self-dual code of length $n$. We also refer to \cite{HuangYNL,ZhangZT} for some recent studies on NMDS self-dual codes. Based on these results, it is interesting and widely open to construct NMDS self-dual codes with new parameters.

\begin{table}[t]\label{table-1}
\caption{Some known results on $q$-ary NMDS self-dual codes}
\centering
\begin{tabular}{|c|c|c|c|c|}
\hline
length $n$ & condition on $q$ & references \\
[0.5ex]
\hline
$n=q-1$& $q\equiv 1\pmod 4$ and $q\le 113$ &\cite{GKL}\\
\hline
$n\le 16$& prime $q\le 197$ &\cite{KotKS}\\
\hline
$n$ with $4|n$ and $q\ge 4^{n+3}(n+3)^2$&$q\equiv 1\pmod 4$&\cite{JinK}\\
\hline
$n$ with $2|n$ and $4\le n\le q+\lfloor 2\sqrt{q}\rfloor-2$& $2|q$ and $q\ge 4$&\cite{JinK}\\
\hline
\end{tabular}
\end{table}

\begin{table}[t]\label{table2}
\caption{Our results on $q$-ary NMDS self-dual codes}
\centering
\begin{tabular}{|c|c|c|c|c|}
\hline
length $n$ & condition on $q$ (odd) & references \\
[0.5ex]
\hline
$ n|(q-1),\ n<q-1$& $q\equiv 1\pmod 4$ & Theorem \ref{multi1} \\
\hline
$n=tf$ and $1\le t\le R$& $q=r^2$, $q-1=ef$ &Theorem \ref{multi2}\\
$R=\frac{r+1}{\gcd(r+1,f)}$&$e$ is even&\\
\hline
$n=s(r-1)+t(r+1)$ &$q=r^2$&Theorem \ref{multi3}\\  $1\le s\le \frac{r+1}{2}$, $1\le t\le \frac{r-1}{2}$ & $r\equiv 1\pmod 4$ and $s$ is even &\\
\hline
$n=s(r-1)+t(r+1)$ &$q=r^2$&Theorem \ref{multi3}\\  $1\le s\le \frac{r+1}{2}$, $1\le t\le \frac{r-1}{2}$ & $r\equiv 3\pmod 4$ and $s$ is odd &\\
\hline
$n=2tr^{\ell}$ &$q=p^m$ and $m$ is even&Theorem \ref{additive1}\\  $0\le\ell<\frac{m}{s}$, $1\le t\le\frac{r-1}{2}$ &  $r=p^s$ with $s|\frac{m}{2}$ &\\
\hline
$n=tr+sp^{t'}$, $t'=\lceil\log_p(t)\rceil$ &$q=r^2$ and $r=p^m$&Theorem \ref{additive2}\\  $1\le t\le r$, $0\le s\le p^{m-t'}-1$ &  $t$ and $s$ even &\\
\hline
\end{tabular}
\end{table}

In this paper, we construct the following class of codes.

\begin{definition} Let $n$ and $k$ be positive integers such that $n>k$.
Let $A=\{a_1,\cdots,a_{n}\}\subset\mathbb F_q$ and $\lambda=(\lambda_1,\lambda_2,\cdots,\lambda_{n})$, where $\lambda_1,\cdots,\lambda_{n}\in\mathbb F_q^*=\mathbb F_q\setminus\{0\}$. We denote by $\mathcal C(A,k,\lambda)$ the $q$-ary linear code with the generator matrix
\begin{equation}\label{Matrix1}
    \begin{pmatrix}
\lambda_1a_1^k & \lambda_2a_2^k & \cdots & \lambda_{n}a_{n}^k\\
\lambda_1a_1^{k-2} & \lambda_2a_2^{k-2} & \cdots & \lambda_{n}a_{n}^{k-2}\\
\vdots & \vdots & \ddots & \vdots\\
\lambda_1a_1^2 & \lambda_2a_2^2 & \cdots & \lambda_{n}a_{n}^2\\
\lambda_1a_1 & \lambda_2a_2 & \cdots & \lambda_{n}a_{n}\\
\lambda_1 & \lambda_2 & \cdots & \lambda_{n}
\end{pmatrix}.
\end{equation}
\end{definition}

We shall prove that $\mathcal C(A,k,\lambda)$ is either MDS or NMDS which depends on certain zero-sum condition. (Proposition \ref{NMDS1}) In the case of NMDS codes, we provide an effective approach (Lemma \ref{NMDSkeylemma}) to construct NMDS self-dual codes which largely extend known parameters of such codes (see Table 2). In particular, we proved that for square $q$, almost $q/8$ NMDS self-dual codes over $\mathbb F_q$ can be constructed.



The following sections are organized as follows. In Section 2, we recall some basic results on MDS codes and NMDS codes. Then we prove that $\mathcal C(A,k,\lambda)$ is either MDS or NMDS which depends on certain zero-sum condition. (Proposition \ref{NMDS1}) In Section 3, we consider the NMDS case of $\mathcal C(A,k,\lambda)$ (Lemma \ref{NMDSkeylemma}) and construct NMDS self-dual codes with new parameters.


\section{On the code $\mathcal C(A,k,\lambda)$}

Firstly, we claim that the code $\mathcal C(A,k,\lambda)$ is an $[n,k]$ code over $\mathbb F_q$, where $A=\{a_1,\cdots,a_{n}\}\subset\mathbb F_q$ and $\lambda=(\lambda_1,\lambda_2,\cdots,\lambda_{n})$ with $\lambda_i\in\mathbb F_q^*$. We have to show that the $k$ rows (denoted by $\alpha_1,\cdots,\alpha_k$) of (\ref{Matrix1}) are linear independent. It suffices to show that the first row $\alpha_1$ is not a linear combination of the rest $k-1$ rows. Assume to the contrary that there exist $c_2,\cdots,c_k\in\mathbb F_q$ such that $\alpha_1=c_2\alpha_2+\cdots+c_k\alpha_k$. This means that $a_1,\cdots,a_n$ are distinct roots of the polynomial
$$f(x)=x^{k}+c_2x^{k-2}+\cdots+c_k.$$
But this is a contradiction, as $n>k$. This proves the claim.

Let $\mathcal C$ be a linear code over $\mathbb F_q$ and $G=(a_{ij})_{k\times n}$ its generator matrix. Let $S=\{\alpha_1,\ldots,\alpha_n\}$ be the multiset of all columns of $G$. Recall the following properties of MDS and NMDS codes.

\begin{proposition}{\rm(\cite{Ball}, Lemma 7.3)}\label{BasicMDS}
The linear code $\mathcal C$ is MDS if and only if any $k$ columns of $G$ are linearly independent.
\end{proposition}

\begin{proposition}{\rm(\cite{DodL})}\label{BasicNMDS}
The linear code $\mathcal C$ is NMDS if and only if the following hold:
\begin{enumerate}
    \item any $k-1$ columns of $G$ are linearly independent;
    \item there exists $k$ columns of $G$ that are linearly dependent;
    \item any $k+1$ columns of $G$ have rank $k$.
\end{enumerate}
\end{proposition}

The following lemma is well known in finite geometry.
\begin{lemma}{\rm(\cite{Ball}, Lemma 5.16)}\label{Hyperplane}
Let $u=(u_1,\ldots,u_k)$ be a non-zero vector of $\mathbb F_q^k$. The hyperplane
$$\ker(u_1x_1+\cdots+u_kx_k)$$
contains $|S|-w$ points of $S$ if and only if the codeword $uG$ has weight $w$.
\end{lemma}

Combining Proposition \ref{BasicNMDS} and Lemma \ref{Hyperplane}, we have the following.

\begin{lemma}\label{NMDSdef}
The linear code $\mathcal C$ is NMDS if and only if the following hold:
\begin{enumerate}
    \item any $k-1$ columns of $G$ are linearly independent;
    \item $|H\cap S|\le k$ for any hyperplane $H$;
    \item there exists a hyperplane $H$ such that $|H\cap S|=k$.
\end{enumerate}
\end{lemma}

Let $A=\{a_1,\ldots,a_n\}\subset\mathbb F_q$. We call $A$
\begin{itemize}
    \item a zero-sum subset of $\mathbb F_q$ if $a_1+\cdots+a_n=0$;
    \item $k$-zero-sum free if $A$ contains no zero-sum subset of size $k$.
\end{itemize}

\begin{proposition}\label{NMDS1}
Let $n$ and $k$ be positive integers such that $n>k$.
Let $A=\{a_1,\cdots,a_{n}\}\subset\mathbb F_q$ and $\textbf{1}=(1,1,\cdots,1)\in \mathbb F_q^n$.
Then $\mathcal C(A,k,\textbf{1})$ is either MDS or NMDS which depends on the following:
\begin{enumerate}
\item $\mathcal C(A,k,\textbf{1})$ is MDS if and only if $A$ is $k$-zero-sum free.
\item $\mathcal C(A,k,\textbf{1})$ is NMDS if and only if $A$ contains a zero-sum subset of size $k$.
\end{enumerate}
\end{proposition}

\begin{proof}(1).
For any $k$ columns $\alpha_{i_1},\cdots,\alpha_{i_{k}}$ of the matrix (\ref{Matrix1}) (where $\lambda=\textbf{1}$), let
$$
B=
\begin{pmatrix}
a_{i_1}^k & a_{i_2}^k & \cdots & a_{i_k}^k\\
a_{i_1}^{k-2} & a_{i_2}^{k-2} & \cdots & a_{i_k}^{k-2}\\
\vdots & \vdots & \ddots & \vdots\\
a_{i_1}^2 & a_{i_2}^2 & \cdots & a_{i_k}^2\\
a_{i_1} & a_{i_2} & \cdots & a_{i_k}\\
1 & 1 & \cdots & 1
\end{pmatrix}.
$$
It is easy to see that $$\det
B=(a_{i_1}+\cdots+a_{i_k})\prod_{1\le s<t\le k}(a_{i_t}-a_{i_s}).$$
By Proposition \ref{BasicMDS}, we have the desired result.

(2). ($\Leftarrow$): We may assume that $a_1+\cdots+a_k=0$. We prove that the matrix (\ref{Matrix1}) (where $\lambda=\textbf{1}$) satisfies (1)-(3) of Lemma \ref{NMDSdef}.

Any $k-1$ columns $\alpha_{i_1},\cdots,\alpha_{i_{k-1}}$ forms a $k\times(k-1)$ matrix which contains an invertible $(k-1)\times(k-1)$ submatrix
$$
\begin{pmatrix}
a_{i_1}^{k-2} & a_{i_2}^{k-2} & \cdots & a_{i_{k-1}}^{k-2}\\
\vdots & \vdots & \ddots & \vdots\\
a_{i_1}^2 & a_{i_2}^2 & \cdots & a_{i_{k-1}}^2\\
a_{i_1} & a_{i_2} & \cdots & a_{i_{k-1}}\\
1 & 1 & \cdots & 1
\end{pmatrix}.
$$
Consequently, $\alpha_{i_1},\cdots,\alpha_{i_{k-1}}$ are linearly independent.

Let
$$H=\ker(u_1x_1+\cdots+u_kx_k)$$
be any hyperplane. Since $f(x)=u_1x^k+u_2x^{k-2}\cdots+u_{k-1}x+u_k$ has at most $k$ distinct roots in $\mathbb F_q$, we have $|H\cap S|\le k$.

For $2\le \ell\le k$, let
$$e_{\ell}=(-1)^{\ell}\sum_{1\le i_1<\cdots<i_{\ell}\le k}a_{i_1}\cdots a_{i_l}$$
and consider the hyperplane
$$H=\ker(x_1+e_2x_2\cdots+e_kx_k).$$
As $a_1+\cdots+a_k=0$, it is easy to see that $|H\cap S|=|\{\alpha_1,\ldots,\alpha_k\}|= k$.

($\Rightarrow$): Since $\mathcal C(A,k,\textbf{1})$ is an NMDS code, there exists $k$ columns $\alpha_{i_1},\cdots,\alpha_{i_{k}}$ of the matrix (\ref{Matrix1}) (where $\lambda=\textbf{1}$) that are linearly dependent. This means the following $k\times k$ submatrix
$$
B=
\begin{pmatrix}
a_{i_1}^k & a_{i_2}^k & \cdots & a_{i_k}^k\\
a_{i_1}^{k-2} & a_{i_2}^{k-2} & \cdots & a_{i_k}^{k-2}\\
\vdots & \vdots & \ddots & \vdots\\
a_{i_1}^2 & a_{i_2}^2 & \cdots & a_{i_k}^2\\
a_{i_1} & a_{i_2} & \cdots & a_{i_k}\\
1 & 1 & \cdots & 1
\end{pmatrix}
$$
is not invertible. Therefore, $$\det
B=(a_{i_1}+\cdots+a_{i_k})\prod_{1\le s<t\le k}(a_{i_t}-a_{i_s})=0,$$
which means $a_{i_1}+\cdots+a_{i_k}=0$.
\end{proof}

\begin{remark}
It is easy to see that the results in Proposition \ref{NMDS1} also hold for the code $\mathcal C(A,k,\lambda)$ for any $\lambda=(\lambda_1,\lambda_2,\cdots,\lambda_{n})$, where $\lambda_1,\cdots,\lambda_{n}\in\mathbb F_q^*$.
\end{remark}

\section{Constructions of NMDS self-dual codes}

Let $\eta(x)$ be the quadratic character of $\mathbb F_q^*$. For any $A\subset \mathbb F_q$ and $a\in A$, we denote
$$\pi_A(a)=\prod_{e\in A,\ e\neq a}(a-e).$$
Recall the following result on MDS self-dual codes.

\begin{lemma}{\rm(\cite{JinX}, Corollary 2.14)}\label{MDS}
Let $n$ be an even number and $A\subset \mathbb F_q$ with $n$ elements.
Assume that $\eta(\pi_A(a))$ are the same for all $a\in A$. Then there exists a $q$-ary MDS self-dual code $\mathcal C$ of length $n$.
\end{lemma}

In literature, most constructions on MDS self-dual codes are based on generalized Reed-Solomon codes. Let $a_1,\cdots,a_n$ be distinct elements in $\mathbb F_q\cup\{\infty\}$ and $\lambda_1,\cdots,\lambda_n\in\mathbb F_q^*$. Denote $A=\{a_1,\cdots,a_n\}$ and $\lambda=(\lambda_1,\cdots,\lambda_n)$. The generalized Reed-Solomon code associated to $A$ and $\lambda$ is defined by
$$GRS_k(A,\lambda)=\{(\lambda_1f(a_1),\ldots,\lambda_nf(a_n))\ |\ f(x)\in\mathbb F_q[x],\ \deg f<k\},$$
where $f(\infty)$ is defined as the coefficient of $x^{k-1}$ in $f$. When $\textbf{1}=(1,\cdots,1)$, the generalized Reed-Solomon code $GRS_k(A,\textbf{1})$ is called a Reed-Solomon code. It is well-known that generalized Reed-Solomon codes are MDS codes.

Based on our construction of the linear code $\mathcal C(A,k,\lambda)$, we have the following crucial lemma, which is an analog of Lemma \ref{MDS} for NMDS self-dual codes.

\begin{lemma}\label{NMDSkeylemma}
Let $A=\{a_1,\cdots,a_{2k}\}\subset\mathbb F_q$ with $a_1+\cdots+a_{2k}=0$.
Assume that $\eta(\pi_A(a))$ are the same for all $a\in A$. Then
\begin{enumerate}
\item if $A$ is $k$-zero-sum free, then there exists $\lambda=(\lambda_1,\cdots,\lambda_n)\in\mathbb F_q^n$ such that $\mathcal C(A,k,\lambda)$ is a $q$-ary MDS self-dual code.
\item if $A$ contains a zero-sum subset of size $k$, then there exists $\lambda=(\lambda_1,\cdots,\lambda_n)\in\mathbb F_q^n$ such that $\mathcal C(A,k,\lambda)$ is a $q$-ary NMDS self-dual code.
\end{enumerate}
\end{lemma}

\begin{proof}
Let $\lambda_1,\lambda_2,\cdots,\lambda_{2k}\in\mathbb F_q^*$ and $\lambda=(\lambda_1,\lambda_2,\cdots,\lambda_{2k})$. It is easy to see that $\mathcal C(A,k,\lambda)$ is self-dual if and only if
$$\sum_{i=1}^n\lambda_i^2a_i^j=0$$
for any $j=0,1,\cdots,2k-2,2k$. That is, $(\lambda_1^2,\lambda_2^2,\cdots,\lambda_{2k}^2)$ are solutions of the following equation
\begin{equation}\label{equation1}
\begin{pmatrix}
a_1^{2k} & a_2^{2k} & \cdots & a_{2k}^{2k}\\
a_1^{2k-2} & a_2^{2k-2} & \cdots & a_{2k}^{2k-2}\\
\vdots & \vdots & \ddots & \vdots\\
a_1^2 & a_2^2 & \cdots & a_{2k}^2\\
a_1 & a_2 & \cdots & a_{2k}\\
1 & 1 & \cdots & 1
\end{pmatrix}
\begin{pmatrix}
x_1\\
x_2\\
\vdots\\
x_{2k-2}\\
x_{2k-1}\\
x_{2k}
\end{pmatrix}=\begin{pmatrix}
0\\
0\\
\vdots\\
0\\
0\\
0
\end{pmatrix}.
\end{equation}

Let
$$B=
\begin{pmatrix}
a_1^{2k} & a_2^{2k} & \cdots & a_{2k}^{2k}\\
a_1^{2k-2} & a_2^{2k-2} & \cdots & a_{2k}^{2k-2}\\
\vdots & \vdots & \ddots & \vdots\\
a_1^2 & a_2^2 & \cdots & a_{2k}^2\\
a_1 & a_2 & \cdots & a_{2k}\\
1 & 1 & \cdots & 1
\end{pmatrix}.$$
Since $a_1+\cdots+a_{2k}=0$, we have $$\det
B=(a_1+\cdots+a_{2k})\prod_{i<j}(a_j-a_i)=0,$$
and consequently (\ref{equation1}) has non-zero solutions. Let $(y_1,\cdots,y_{2k})$ be a non-zero solution of (\ref{equation1}), then we have $y_i\neq 0$ for all $1\le i\le 2k$. We may assume that $y_{2k}=1$. As $\text{Rank}(B)=2k-1$, we get the following solution
\begin{equation}\label{solution1}
    y_1=-\frac{\prod_{j=2}^{2k-1}(a_j-a_{2k})}{\prod_{j=2}^{2k-1}(a_j-a_1)},
\cdots,
y_{2k-1}=-\frac{\prod_{j=1}^{2k-2}(a_j-a_{2k})}{\prod_{j=1}^{2k-2}(a_j-a_{2k-1})}.
\end{equation}
From (\ref{solution1}), we obtain the following solution
$$y_1=\prod_{j=2}^{2k}(a_1-a_j)^{-1}=\pi_A(a_1)^{-1},
\cdots,y_{2k}=\prod_{j=1}^{2k-1}(a_{2k}-a_j)^{-1}=\pi_A(a_{2k})^{-1}.$$
As $\eta(\pi_A(a))$ are the same for all $a\in A$, we can obtain
$(\lambda_1,\lambda_2,\cdots,\lambda_{2k})$ from $(y_1,\ldots,y_{2k})$. Therefore, $\mathcal C(A,k,\lambda)$ is a self-dual code.

(1). In this case, it follows easily from Proposition \ref{NMDS1}.(1) that $\mathcal C(A,k,\lambda)$ is a $q$-ary MDS self-dual code.

(2). In this case, it follows easily from Proposition \ref{NMDS1}.(2) that $\mathcal C(A,k,\lambda)$ is a $q$-ary NMDS self-dual code.
\end{proof}

In the following, with the help of Lemma \ref{NMDSkeylemma}, combining with some known constructions on MDS self-dual codes, we construct more NMDS self-dual codes. Roughly speaking, most constructions of MDS self-dual codes are based on unions of cosets of multiplicative subgroups or additive subgroups. Let $QR_q$ be the set of all nonzero squares of $\mathbb F_q$.

\subsection{Constructions via multiplicative subgroups}
\ \medskip

In the following, we choose some representative constructions \cite{FLL,HYan,ZhangFeng}, which are constructed using multiplicative subgroups.

\subsubsection{Constructions via one multiplicative subgroups}

\begin{theorem}\label{multi1}
Let $q$ be an odd prime power and $n$ an even positive integer. Then there exists a $q$-ary NMDS self-dual code of length $n$ if $q$ and $n$ satisfy the following condition:
$$ n|(q-1),\ n<q-1,\text{ and }q\equiv 1\pmod 4.
$$
\end{theorem}

\begin{proof}
{\bf Case 1:} Assume that $n\equiv 2\pmod 4$. Let $\theta$ be a primitive $n$-th root of unity of $\mathbb F_q^*$. Let $S=\{\theta,\theta^2,\cdots,\theta^{n}\}$. It is verified in \cite[Theorem 1]{HYan} that $\eta(\pi_S(a))$ are the same for all $a\in S$. Let $n=2m$.
Since
$$x^{2m}-1=\prod_{i=1}^{2m}(x-\theta^{i})
\text{ and }x^m-1=\prod_{i=1}^m(x-\theta^{2i}),$$
we have
$\sum_{i=1}^{2m}\theta^i=0 \text{ and } \sum_{i=1}^{m}\theta^{2i}=0$.
By Lemma \ref{NMDSkeylemma}, we obtain the desired result.

{\bf Case 2:} Assume that $n\equiv 0\pmod 4$ and $n=2m$. Let $\theta$ be a primitive $m$-th root of unity of $\mathbb F_q^*$. Choose $\beta\in QR_q\setminus\{\theta,\theta^2,\cdots,\theta^{m}\}$. Let $$S=\{\theta,\theta^2,\cdots,\theta^{m},\beta\theta,\beta\theta^2,\cdots,\beta\theta^m\}.$$ It is verified in \cite[Theorem 1]{HYan} that $\eta(\pi_S(a))$ are the same for all $a\in S$. Since
$$(x^{m}-1)(x^m-\beta^m)=\prod_{i=1}^{m}(x-\theta^{i})\prod_{i=1}^{m}(x-\beta\theta^{i}),$$
we have
$\sum_{a\in S}a=0 \text{ and } \sum_{i=1}^{m}\theta^{i}=0$.
By Lemma \ref{NMDSkeylemma}, we obtain the desired result.
\end{proof}

\begin{theorem}\label{multi2}
Let $r$ be an odd prime power, $q=r^2$, $q-1=ef$ and $R=\frac{r+1}{\gcd(r+1,f)}$.
Assume that $tf$ is even and $1\le t\le R$. If $e$ is even, then there exists a $q$-ary NMDS self-dual code of length $tf$.
\end{theorem}

\begin{proof}
Let $\mathbb F_q^*=\langle\theta\rangle$, $C=\langle\alpha\rangle$, $M=\langle\beta\rangle$ with $\alpha=\theta^e$ and $\beta=\theta^{r-1}$. Assume that $\{i_1,i_2,\ldots,i_t\}$ is a subset of $\mathbb Z_{r+1}$ such that $i_1,i_2,\ldots,i_t$ are distinct module $\frac{r+1}{\gcd(r+1,f)}$. Let $B=\{\beta^{i_{\lambda}}\ :\ 1\le \lambda\le t\}$ and $S=\cup_{\lambda=1}^t\beta^{i_{\lambda}}C$. It is verified in \cite[Theorem 7]{ZhangFeng} that $\eta(\pi_S(a))$ are the same for all $a\in S$. For any $1\le\lambda\le t$, we have
$$\prod_{a\in \beta^{i_{\lambda}}C}(x-a)=x^f-\beta^{i_{\lambda}f}.$$ Therefore, $\sum_{a\in \beta^{i_{\lambda}}C}a=0$.

If $t$ is even, then we have $\sum_{a\in S}a=0$ and $\sum_{\lambda=1}^{t/2}\sum_{a\in \beta^{i_{\lambda}}C}a=0$.

Assume that $t$ is odd and $f$ is even. Let $C'=\langle\alpha^2\rangle$. Then $\cup_{\lambda=1}^t\beta^{i_{\lambda}}C'\subset \cup_{\lambda=1}^t\beta^{i_{\lambda}}C$
and $|\cup_{\lambda=1}^t\beta^{i_{\lambda}}C'|=tf/2$. For any $1\le\lambda\le t$, we have
$$\prod_{a\in \beta^{i_{\lambda}}C'}(x-a)=x^{f/2}-
\beta^{i_{\lambda}f/2}.$$
Therefore, $\sum_{a\in \beta^{i_{\lambda}}C'}a=0$ and $\sum_{\lambda=1}^{t}\sum_{a\in \beta^{i_{\lambda}}C'}a=0$. In all cases, by Lemma \ref{NMDSkeylemma}, we obtain the desired result.
\end{proof}

\subsubsection{Constructions via two multiplicative subgroups}
\ \medskip

An important and striking construction is due to Fang, Liu and Luo \cite{FLL}, who constructed, for large square $q$, almost $\frac{1}{8}q$ $q$-ary MDS self-dual codes. Now we employ their construction and verify the conditions in Lemma \ref{NMDSkeylemma}. Actually, we provide a slightly simplified proof. For any two elements $a,b\in\mathbb F_q^*$, we define
$a\approx b$ if there exists $c\in QR_q$ such that $a=bc$.

\begin{theorem}\label{multi3}
Let $q=r^2$, where $r$ is an odd prime power. For any $1\le s\le \frac{r+1}{2}$ and $1\le t\le \frac{r-1}{2}$, assume that $n=s(r-1)+t(r+1)$. There exists a $q$-ary NMDS self-dual code of length $n$ if $r$ and $s$ satisfy one of the following:
\begin{enumerate}
    \item $r\equiv 1\pmod 4$ and $s$ is even;
    \item $r\equiv 3\pmod 4$ and $s$ is odd.
\end{enumerate}
\end{theorem}

\begin{proof}

Let $\mathbb F_{r^2}^*=\langle\alpha\rangle$. Then $\mathbb F_r^*=\langle\alpha^{r+1}\rangle$. We denote $\gamma=\alpha^{r+1}$ and $\beta=\alpha^{r-1}$.

Consider the following subset
$$A=\alpha^2\langle \gamma\rangle\cup \alpha^4\langle \gamma\rangle \cup\cdots\cup \alpha^{2s}\langle \gamma\rangle\cup\alpha\langle \beta\rangle\cup \alpha^3\langle \beta\rangle \cup\cdots\cup \alpha^{2t-1}\langle \beta\rangle,$$
where $1\le s\le\frac{r+1}{2}$ and $1\le t\le \frac{r-1}{2}$.

It is easy to see that $\alpha^{2i}\langle \gamma\rangle \cap \alpha^{2j}\langle \gamma\rangle =\emptyset$ for $1\le i<j\le\frac{r+1}{2}$ and $\alpha^{2i-1}\langle \beta\rangle \cap \alpha^{2j-1}\langle \beta\rangle=\emptyset$ for $1\le i<j\le \frac{r-1}{2}$. As $\gamma\in QR_q$ and $\alpha,\beta\notin QR_q$, we have $\alpha^{2i}\langle \gamma\rangle \cap \alpha^{2j-1}\langle \beta\rangle=\emptyset$. In the following, we show that $\eta(\pi_A(a))$ are the same for any $a\in A$.

For any $1\le i\le\frac{r+1}{2}$ and $1\le j\le \frac{r-1}{2}$, we denote
$$f_i(x):=\prod_{a\in \alpha^{2i}\langle \gamma\rangle}(x-a)\text{ and }
g_j(x):=\prod_{a\in \alpha^{2j-1}\langle \beta\rangle}(x-a).$$
Note that
\begin{equation}\label{Luofg}
    f_i(x)=x^{r-1}-\alpha^{2i(r-1)}\text{ and } g_j(x)=x^{r+1}-\alpha^{(2j-1)(r+1)}.
\end{equation}
Let $a\in\alpha^{2i}\langle \gamma\rangle$. It is easy to see
$$\prod_{e\in \alpha^{2i}\langle \gamma\rangle,\ e\neq a}(a-e)
=f_{i}'(a)=-a^{r-2}\approx 1.$$
For $1\le j\neq i\le\frac{r+1}{2}$, we have
$$\prod_{e\in \alpha^{2j}\langle \gamma\rangle}(a-e)
=f_{j}(a)=a^{r-1}-\alpha^{2j(r-1)}=\alpha^{2i(r-1)}-\alpha^{2j(r-1)}.$$
Note that
$$\alpha^{2i(r-1)}-\alpha^{2j(r-1)}\approx 1-\alpha^{2(j-i)(r-1)}
=1-\beta^{2(j-i)}.$$
We claim that, for any $i\in[-\frac{r+1}{2},\frac{r+1}{2}]\setminus\{0\}$, we have $1-\beta^{2i}\approx \alpha^{\frac{r+1}{2}}$. As $\beta=\alpha^{r-1}$ and $-1=\alpha^{\frac{r-1}{2}(r+1)}$, we have
$$(1-\beta^{2i})^r=1-\beta^{2ir}=1-\beta^{-2i}=\beta^{-2i}(\beta^{2i}-1),$$
and $(1-\beta^{2i})^{r-1}=-\beta^{-2i}=(\alpha^{\frac{r+1}{2}-2i})^{r-1}$. Therefore, we have $\frac{1-\beta^{2i}}{\alpha^{\frac{r+1}{2}-2i}}\in\mathbb F_r$ and the claim follows.

By the above claim, for $1\le j\neq i\le\frac{r+1}{2}$, we have
$$\prod_{e\in \alpha^{2j}\langle \gamma\rangle}(a-e)
\approx \alpha^{\frac{r+1}{2}}.$$

For $1\le j\le \frac{r-1}{2}$, it is easy to see that
$$\prod_{e\in \alpha^{2j-1}\langle \beta\rangle}(a-e)
=a^{r+1}-\alpha^{(2j-1)(r+1)}\in\mathbb F_r.$$

As a result, for any $a\in\alpha^{2i}\langle \gamma\rangle$, we have $$\prod_{e\in A}(a-e)\approx\alpha^{\frac{r+1}{2}(s-1)}.$$

Next, we consider $b=\alpha^{2j-1}\beta^k\in\alpha^{2j-1}\langle \beta\rangle$. Firstly, we have
$$\prod_{e\in \alpha^{2j-1}\langle \beta\rangle,\ e\neq b}(b-e)
=g_{j}'(b)=b^{r}\approx \alpha.$$
For $1\le i\neq j\le \frac{r-1}{2}$, similar to the above, we have
$$\prod_{e\in \alpha^{2j-1}\langle \beta\rangle}(b-e)
\in\mathbb F_r.$$
For $1\le i\le\frac{r+1}{2}$, we have
\begin{align*}
\prod_{e\in \alpha^{2i}\langle \gamma\rangle}(b-e)
=f_{i}(b)&=b^{r-1}-\alpha^{2i(r-1)}\\
&=\beta^{2j-1}\beta^{k(r-1)}-\beta^{2i}\\
&=\beta^{2j-2k-1}-\beta^{2i}\approx \beta^{2\ell-1}-1
\end{align*}
for some integer $\ell$. Similar to the above, we have
$$\beta^{2\ell-1}-1\approx \alpha^{\frac{r-1}{2}}.$$
Consequently, for any $b\in\alpha^{2j-1}\langle \beta\rangle$, we have $$\prod_{e\in A}(b-e)\approx\alpha^{\frac{r-1}{2}s}\alpha.$$
In order to show that $\eta(\pi_A(a))$ are the same for any $a\in A$, it suffices to verify that $$\frac{r+1}{2}(s-1)\equiv\frac{r-1}{2}s+1\pmod 2,$$
which is implied by condition (1) or (2).

From (\ref{Luofg}), it is easy to see that $\sum_{a\in A}a=0$.
For $r\ge 5$, let
$$A'=\alpha^2\langle \gamma^2\rangle\cup \alpha^4\langle \gamma^2\rangle \cup\cdots\cup \alpha^{2s}\langle \gamma^2\rangle\cup\alpha\langle \beta^2\rangle\cup \alpha^3\langle \beta^2\rangle \cup\cdots\cup \alpha^{2t-1}\langle \beta^2\rangle,$$
where $1\le s\le\frac{r+1}{2}$ and $1\le t\le \frac{r-1}{2}$. Then $A'\subset A$ and $|A'|=|A|/2$.
For any $1\le i\le\frac{r+1}{2}$ and $1\le j\le \frac{r-1}{2}$, we have
$$\prod_{a\in \alpha^{2i}\langle \gamma^2\rangle}(x-a)
=x^{\frac{r-1}{2}}-\alpha^{2i\frac{r-1}{2}}
\text{ and }
\prod_{a\in \alpha^{2j-1}\langle \beta\rangle}(x-a)
=x^{\frac{r+1}{2}}-\alpha^{(2j-1)\frac{r+1}{2}}.$$
Consequently, we have $\sum_{a\in A'}a=0$. The remaining single case is when $r=3$, i.e., $q=9$.
In this case, $A=\{\alpha,\alpha^{2},\alpha^{3},\alpha^{5},\alpha^{6},\alpha^{7}\}$. As $\alpha^4=-1$, we may assume that $\alpha$ satisfy $\alpha^2=\alpha+1$. Then we choose $A'=\{\alpha,\alpha^2,\alpha^7\}$. Therefore $|A'|=|A|/2$ and $\sum_{a\in A'}a=0$. In both cases, by Lemma \ref{NMDSkeylemma}, we obtain the desired result.
\end{proof}

\subsection{Constructions via additive subgroups}
\ \medskip

In the following, we choose some other representative constructions \cite{FangF,FangXF}, which are constructed using additive subgroups.

\begin{theorem}\label{additive1}
Let $q=p^m$ be an odd prime power. Assume that $m$ is even and $r=p^s$ with $s|\frac{m}{2}$. Let $n=2tr^{\ell}$, where $0\le\ell<\frac{m}{s}$ and $1\le t\le\frac{r-1}{2}$. Then there exists a $q$-ary NMDS self-dual code of length $n$.
\end{theorem}

\begin{proof}
Let $H$ be an $\mathbb F_r$-subspace of $\mathbb F_q$ of dimension $\ell$ and let $\alpha\in \mathbb F_q\setminus H$. Let $\mathbb F_r=\{\xi_0=0,\xi_1=1,\ldots,\xi_{r-1}\}$ and $H_i=H+\xi_i\alpha$ for $0\le i\le r-1$. Let $S=\cup_{i=0}^{2t-1}H_i=\{a_1,a_2,\ldots,a_n\}$. It is verified in \cite[Theorem 3.3]{FangF} that $\eta(\pi_S(a))$ are the same for all $a\in S$. It is easy to see that $\sum_{a\in S}a=0$ and $\sum_{a\in \cup_{i=0}^{t-1}H_i}a=0$.
By Lemma \ref{NMDSkeylemma}, we obtain the desired result.
\end{proof}

\begin{theorem}\label{additive2}
Let $r=p^m$ and $q=r^2$, where $p$ is an odd prime. For any even $t$ and even $s$ with $1\le t\le r$ and $0\le s\le p^{m-t'}-1$, let $n=tr+sp^{t'}$, where $t'=\lceil\log_p(t)\rceil$. Then there exists a $q$-ary NMDS self-dual code of length $n$.
\end{theorem}

\begin{proof}
Let $H$ be an $\mathbb F_p$-linear subspace of $\mathbb F_r$ of dimension $t'$. Then $|H|=p^{t'}\ge t$ and $|\mathbb F_r/H|=p^{m-t'}>s$. Let $h_1=0,h_2,\ldots,h_t$ be $t$ distinct elements of $H$. Let $b_0=0,b_1,\ldots,b_s$ be $s+1$ distinct representations of $\mathbb F_r/H$ such that for any $1\le i\le \frac{s}{2}$, $b_i=-b_{\frac{s}{2}+i}$. For any $1\le i\le t$, define $T_i=\{x\in\mathbb F_q\ :\ \operatorname{Tr}(x)=h_i\}$. Then $|T_i|=r$ and $T_i\cap T_j=\emptyset$ for any $i\neq j$. For any $0\le j\le s$, define $H_j=\{b_j+h\ :\ h\in H\}$. Let $S=(\cup_{i=1}^tT_i)\cup(\cup_{j=1}^sH_j)$. It is verified in \cite[Theorem 1]{FangXF} that $\eta(\pi_S(a))$ are the same for all $a\in S$. It is easy to see that $\sum_{a\in S}a=0$. Since $\prod_{a\in T_i}(x-a)=x^r+x-h_i$, we have $\sum_{a\in T_i}a=0$. Let $S'=(\cup_{i=1}^{t/2}T_i)\cup(\cup_{j=1}^{s/2}H_j)$. Then we have $|S'|=|S|/2$ and $\sum_{a\in S'}a=0$. By Lemma \ref{NMDSkeylemma}, we obtain the desired result.
\end{proof}

\begin{table}[t]\label{table-3}
\caption{Proportion of number of possible lengths $N$ to $q=r^2$}
\centering
\begin{tabular}{|c|c|c|c|c|}
\hline
$r$ & $q$ & $N$ & $N/q$ \\
[0.5ex]
\hline
101 & 10201 & 1528 & $14.97\%$\\
\hline
107 & 11449 & 1586 & $13.85\%$\\
\hline
199 & 39601 & 5211 & $13.15\%$\\
\hline
\end{tabular}
\end{table}

As we have mentioned that, the maximal length for $q$-ary nontrivial NMDS codes is conjectured to be at most $2q + 2$. It is easy to see that (by Theorem \ref{multi3} and \cite[Remark III.1]{FLL}), our constructions lead to almost $\frac{q/8}{2q/2}=\frac{1}{8}$ of all possible lengths of $q$-ary NMDS self-dual code, which largely extend known results \cite{JinK} on such codes. In Table 3, we list some examples.






\subsection*{Acknowledgments}

D. Han's research was supported by the National Natural Science Foundation of China under Grant 12171398, the Natural Science Foundation of Sichuan Province under Grant 2022NSFSC1856 and the Fundamental Research Funds for the Central Universities under Grant 2682020ZT101. H.B. Zhang was supported by the National Science Foundation of China Grant No. 11901563 and Guangdong Basic and Applied Basic Research Foundation Grant No. 2021A1515010216.

\end{document}